\documentclass[a4paper, 11pt]{article}
\usepackage[utf8]{inputenc}
\usepackage{color}
\usepackage{amssymb, amsmath, mathbbol, amsfonts}
\usepackage{latexsym}
\usepackage{mathrsfs}
\usepackage{bbm}
\usepackage{amsthm}
\usepackage[left=2cm, right=2cm, top=3cm, bottom=3cm]{geometry}
\usepackage{titleref}
\usepackage{prettyref}
\usepackage{graphicx}
\usepackage{enumerate}
\usepackage{trfsigns}
\usepackage{tikz}
\usepackage{tabularx}
\usepackage{epsfig}
\usepackage{enumerate}
\usepackage[active]{srcltx}
\usetikzlibrary{decorations,arrows}
\usetikzlibrary{decorations.pathmorphing}
\usepgflibrary{decorations.pathreplacing}

\newcommand{\h}{H(K_{s,{\alpha},{\boldsymbol{\gamma}}})}

\newcommand{\boldzero}{\boldsymbol{0}}
\newcommand{\bsgamma}{\boldsymbol{\gamma}}

\newcommand{\boldh}{\boldsymbol{h}}

\newcommand{\boldp}{\boldsymbol{p}}

\newcommand{\bolds}{\boldsymbol{s}}

\newcommand{\boldx}{\boldsymbol{x}}
\newcommand{\boldy}{\boldsymbol{y}}
\newcommand{\boldz}{\boldsymbol{z}}

\newcommand{\NN}{\mathbbm{N}}
\newcommand{\ZZ}{\mathbbm{Z}}

\newcommand{\ii}{\mathrm{i}}

\newcommand{\znw}[1]{\mathcal{Z}_{N,w_{#1{}}}}
\newcommand{\zn}{\mathcal{Z}_{N}}

\newcommand{\ewortd}{e^2(\tilde{z}_1, \dots , \tilde{z}_d)}

\newcommand{\ewortdz}{e^2(\tilde{z}_1, \dots , \tilde{z}_d, z)}
\newcommand{\ewordd}{e^2(z_1, \dots , z_d, z_{d+1})}

\newcommand{\eworts}{e^2(\tilde{z}_1, \dots , \tilde{z}_s)}
\newcommand{\ewortdd}{e^2(\tilde{z}_1, \dots , \tilde{z}_d , \tilde{z}_{d+1})}

\newtheorem{thm}{Theorem}[section]

\newtheorem{alg}{Algorithm}
\newtheorem{rem}{Remark}

\newenvironment{proof1}{\begin{trivlist}\item[\hskip\labelsep{\it Proof.}]}{$\hfill\Box$\end{trivlist}}

\allowdisplaybreaks

\title{On combined component-by-component constructions of lattice point sets}
\author{Helene Laimer\thanks{H. Laimer is supported by the Austrian 
Science Fund (FWF): Projects F5506-N26 and F5508-N26, which are part of the 
Special Research Program "Quasi-Monte Carlo Methods: Theory and Applications".}}
\date{Dedicated to the memory of Joseph F. Traub (1932-2015)}

\begin{document}
\setcounter{page}{1}
\newpage
\pagestyle{plain}
\maketitle

\begin{abstract}
\noindent
The standard method for constructing generating vectors for good lattice point sets is the component-by-component construction. Numerical experiments have shown that the generating vectors found by these constructions sometimes tend to have recurring components, which can lead to the problem of having projections with all lattice points lying on the main diagonal. In this paper we combine methods of Dick and Kritzer to avoid this problem with a reduced fast component-by-component construction. That is, we give a variation of the standard component-by-component construction which avoids repeated components and simultaneously results in a considerable speed-up in comparison to the standard construction.
\end{abstract}

\section{Introduction}

A rank-1 lattice point set consisting of $N$ points is usually constructed with the aid of a generating vector $\boldz$ as follows. Let $\boldz = (z_1, \dots, z_s) \in \{ 1, \dots, N-1 \}^s$. The corresponding rank-1 lattice point set is then given by
\begin{align*}
P = \left\{ \left\{\frac{k \boldz}{N} \right\} \colon k = 0, \dots, N-1 \right\} = \left\{ \boldp_0, \dots, \boldp_{N-1} \right\},
\end{align*}
where the braces around $k \boldz/N$ indicate that we consider the fractional part of each coordinate of $k \boldz/N$. For dimensions greater than 2, however, we do not have an explicit construction for generating vectors which produce in some sense good lattice point sets e.g., point sets that can be used as the integration nodes in quasi-Monte Carlo algorithms. Therefore one usually resorts to computer search algorithms, in particular to so-called component-by-component (for short CBC) algorithms, which choose one component of $\boldz$ at a time. For detailed information see for example \cite{DKS,J,NC2}. To ensure that lattice points, constructed in this manner, are evenly spread, it is common to use only generating vectors $\boldz$ whose components and $N$ are co-prime, i.e., each component of $\boldz$ is chosen from 
\begin{align*}
\mathcal{Z}_{N} := \left\{ z \in \{ 1, \dots , N-1 \} \colon \gcd(z,N) = 1 \right\}.
\end{align*}

In the present paper, however, we want to reduce the search space even further. As in \cite{DKLP} we want to adjust the size of the search space for each coordinate according to its importance. This is called a reduced CBC construction.

Let $s \in \mathbbm{N}$ be the dimension, let $N = b^m$ be the number of points in our lattice point set, with $m \in \mathbbm{N}$ and $b \in \mathbbm{P}$, where $\mathbbm{P}$ denotes the set of all primes. Further let $0 = w_1 \leq w_2 \leq w_3 \leq \cdots$, with $w_j \in \NN_{0}$ for all $j \in \NN$, be the sequence which models the importance of the coordinates of our generating vector. The reduced search spaces are then defined as
\begin{align*}
\znw{j} := \begin{cases}
													\left\{ z \in \{ 1, \dots , b^{m-w_j}-1 \} \colon \gcd(z,b^m) = 1 \right\}	& \text{ if } w_j < m,\\
													\{ 1 \} 																																		& \text{ otherwise},
												\end{cases}
\end{align*}
and the $j$-th component of the generating vector $\boldz$ is chosen from the set $\znw{j}$.

This can lead to a considerable speed-up of the run time of the algorithm, which is illustrated by numerical results in \cite{DKLP}. The algorithm can be sped up even further by the so-called fast or reduced fast CBC algorithm. For detailed information see \cite{DKLP, LP, NC1, NC2}.\newline

For the standard as well as the reduced (fast) CBC construction (see \cite{DKS,J} and \cite{DKLP,NC1,NC2}) it is known from \cite{K} that for large dimensions the components of the produced generating vector tend to recur, which can lead to the problem of having projections where all lattice points lie on the main diagonal. Gantner and Schwab \cite{GS} as well as Dick and Kritzer \cite{DK} have come up with methods to avoid this problem. Gantner and Schwab call their method pruning in the CBC algorithm, while Dick and Kritzer name their refined version of this method projection-corrected CBC construction. The general idea is in each step of the CBC algorithm to define some exclusion set $\mathcal{E}$ whose elements cannot be selected in this step. By defining the exclusion sets as the sets consisting of all elements chosen in the previous steps one can effectively avoid components showing up several times. This method can also be used to avoid other phenomena, like for example all lattice points lying on an antidiagonal. For detailed information see \cite{DK}.

In this paper we combine the reduced and the reduced fast, respectively, with the projection-corrected CBC construction to get a construction which pools the advantages of these two constructions, that is, being considerably faster than the standard CBC construction and being free of recurring components.\newline

This paper is organized as follows. In the next section we define the function space whose elements we would like to numerically integrate by quasi-Monte Carlo quadratures, using as integration nodes the lattice points constructed by our algorithm. Finally, in Section~\ref{sec:main} we state our combined CBC algorithm and prove upper bounds on the worst-case error of the resulting generating vectors.

\section{Definition of the function space}
As mentioned before, with our combined CBC algorithm we would like to construct generating vectors for lattice point sets used in quasi-Monte Carlo (or for short QMC) algorithms applied to functions in certain weighted Korobov spaces.

To define these spaces let us first introduce some notation.
For $s \in \NN$ let $[s]$ denote the set $\{ 1, \dots, s \}$. For some vector $\boldx = (x_1, \dots, x_s) \in [0,1]^s$ and some $\mathfrak{u} \subseteq [s]$ we write $\boldx_{\mathfrak{u}} = (x_j)_{j \in \mathfrak{u}} \in [0,1]^{\left \lvert \mathfrak{u} \right \rvert}$ and $(\boldx_{ \mathfrak{u}}, \boldzero )$ for the vector $(\tilde{x}_1, \dots, \tilde{x}_s) \in \ZZ^s$, with
\begin{align*}
\tilde{x}_j = \begin{cases}
								x_j	& \text{ if } j \in \mathfrak{u},\\
								0		& \text{ otherwise.}
							\end{cases}
\end{align*}

For now let us restrict ourselves to product weights $\bsgamma = \left( \gamma_j \right)_{j \in \NN_0}$, with $\gamma_{\mathfrak{u}} = \prod_{j \in \mathfrak{u}} {\gamma_j}$ for any $\mathfrak{u} \subseteq [s]$ to model the importance of coordinates or groups of coordinates for functions from the Korobov space.

Now we are ready to define the weighted Korobov spaces $\h$. These are reproducing kernel Hilbert spaces of functions defined on $[0,1]^s$, with their reproducing kernel given by
\begin{align*}
K_{s, \alpha, \bsgamma}(\boldx, \boldy) = 1 + \sum_{\emptyset \neq \mathfrak{u} \subseteq [s]} \sum_{\boldh_{\mathfrak{u}} \in \ZZ^{\left \lvert \mathfrak{u} \right \rvert} \setminus \{ \boldzero \} } {r_{\alpha}{(\bsgamma, (\boldh_{\mathfrak{u}}, \boldzero))} \exp{(2 \pi \ii \boldh_{\mathfrak{u}} \cdot (\boldx_{\mathfrak{u}} - \boldy_{\mathfrak{u}} ) ) } } ,\quad \boldx, \boldy \in [0,1)^s,
\end{align*}
where ``$\cdot$'' denotes the usual dot product and where, for $\emptyset \neq \mathfrak{u} \subseteq [s]$ and $\boldh_{\mathfrak{u}} \in \ZZ^{\left \lvert \mathfrak{u} \right \rvert} \setminus \{ \boldzero \}$, we have

\begin{align}\label{eq:weight}
r_{\alpha}(\bsgamma, (\boldh_{\mathfrak{u}}, \boldzero)) = \prod_{j \in \,\mathfrak{u}}{ r_{\alpha}(\gamma_j, h_j) },
\end{align}
with
\begin{align}\label{eq:weight1dim}
r_{\alpha}(\gamma_j, h_j) = \begin{cases}
															1 & \text{ if } h_j = 0,\\
															\frac{\gamma_j}{\left \lvert h_j \right \rvert^{\alpha} } & \text{ otherwise}.
														\end{cases}
\end{align}
For $f, g \in \h$ the inner product is then given by
\begin{align*}
\left\langle f, g \right\rangle_{\h} = 1 + \sum_{\emptyset \neq \mathfrak{u} \subseteq [s]} \sum_{\boldh_{\mathfrak{u}} \in \ZZ^{\left \lvert \mathfrak{u} \right \rvert} \setminus \{ \boldzero \} }\left(r_{\alpha}{(\bsgamma, (\boldh_{\mathfrak{u}}, \boldzero))}\right)^{-1}\hat{f}((\boldh_{\mathfrak{u}}, \boldzero))\overline{\hat{g}((\boldh_{\mathfrak{u}}, \boldzero))},
\end{align*}
where $\hat{f}(\boldh) = \int_{[0,1]^{\bolds}}f(\boldx) \exp(-2 \pi \ii \boldh \cdot \boldx) \,\mathrm{d}\boldx$ denotes the $\boldh$-th Fourier coefficient of $f$. The norm in $\h$ is the norm induced by this inner product.

As a quality measure for a generating vector constructed with our algorithm we want to consider the (squared) worst-case error of integration in $\h$ by a QMC rule using the lattice point set as integration nodes. For a generating vector $\boldz = (z_1, \dots, z_s)$ we define the worst-case error of $\boldz$ as
\begin{align*}
e(\boldz) = e(z_1, \dots , z_s) = \sup_{\substack{f \in \h \\ \left\lVert f \right\rVert_{\h} \leq 1}}{ \left\lvert \int_{[0,1]^s}{f(\boldx) \, \mathrm{d}\boldx} - \frac{1}{N} \sum_{j = 0}^{N-1}{f(\boldp_j)} \right\rvert}.
\end{align*}
It is known (see for example \cite{DKLP,DSWW}) that for a generating vector $\boldz \in \{ 0, \dots N-1 \}^s$ the squared worst-case error in the weighted Korobov space $\h$ is given by
\begin{align}\label{eq:squaredWorst}
e^2(\boldz) = \sum_{\substack{\boldh \in \mathbbm{Z}^s\setminus\{ \boldzero \} \\ \boldz \cdot \boldh \equiv 0 \pmod{N}}}{r_{\alpha}(\bsgamma,\boldh)}.
\end{align}

\section{The combined CBC algorithm}\label{sec:main}

Before we describe the combined CBC construction let us first introduce a little more notation. We denote by $t_1 = \max \{j \in \NN \colon w_j = 0 \}$ the index up to which we consider the whole set $\mathcal{Z}_N$ as the search space and by $t_2 = \min \{j \in \NN \colon w_j \geq m \}$ the first index for which the search space reduces to $\{ 1 \}$. Note that $t_1 \geq 1$, $t_2 \geq 2$, and $t_1 < t_2$. Furthermore, let for sets $\mathcal{E} \subseteq \mathbbm{Z}$, $\left\lvert \mathcal{E} \right\rvert$ denote their cardinality.

Now we are ready to state the combined CBC algorithm. 
\begin{alg}\label{alg}
Let $s \in \NN$, $b \in \mathbbm{P}$, $m \in \NN$, $N = b^m$, $0 = w_1 \leq w_2 \leq \cdots, $ and $\mathcal{Z}_N, \znw{j}, t_1$ and $t_2$ as above.
\begin{enumerate}
	\item Set $z_1 = 1$ and set $\mathcal{E}_1 = \emptyset$.
	\item For $d \in \{ 1, \dots , \min\{t_1 - 1,s-1\}\}$ do the following: Assume that $z_1, \dots , z_d$ have already been found and choose $\mathcal{E}_{d+1} \subsetneq \zn$. (If no coordinates are to be excluded in this step, we define $\mathcal{E}_{d+1} = \emptyset$.) Now choose $z_{d+1} \in \zn \setminus \mathcal{E}_{d+1}$ such that
	\begin{align*}
	\ewordd
	\end{align*}
is minimized as a function of $z_{d+1}$. Set $\tilde{z}_{d+1} = z_{d+1}$.
	\item Increase $d$ by 1 and repeat Step 2 until $d=\min\{t_1-1,s-1\}$.
	\item If $t_1 \geq s$ the algorithm terminates with $d+1=s$. Else, for $d \in \{ t_1, \dots , \min\{ t_2 - 2, s-1\}\}$ do the following: Assume that $z_1, \dots , z_{t_1}, z_{t_1+1}, \dots , z_d$ have already been found and choose $\mathcal{E}_{d+1} \subsetneq \znw{d+1}$. (If no coordinates are to be excluded in this step, we define $\mathcal{E}_{d+1} = \emptyset$.) Now choose $z_{d+1} \in \znw{d+1} \setminus \mathcal{E}_{d+1}$ such that
	\begin{align*}
	e^2(z_1, \dots , z_{t_1}, b^{w_{t_{1}+1}}z_{t_1+1}, \dots, b^{w_d}z_d , b^{w_{d+1}}z_{d+1})
	\end{align*}
is minimized as a function of $z_{d+1}$. Set $\tilde{z}_{d+1} = b^{w_{d+1}} z_{d+1}$.
	\item Increase $d$ by 1 and repeat Step 4 until $d = \min\{ t_2 - 2, s-1\}$.
	\item If $t_2 > s$ the algorithm terminates with $d+1=s$. Else, for $d \in \{ t_2 - 1, \dots, s-1\}$ set $z_d = 0.$ (The corresponding exclusion set is the empty set.) Set $\tilde{z}_d = z_d.$
	\item Increase $d$ by 1 and repeat Step 6 until $d = s-1$.
\end{enumerate}
\end{alg}

To avoid lengthy case analyses let us here and, if not stated otherwise, for the rest of this paper, assume that $t_2 \leq s$. The proofs of \prettyref{thm:betterbound} for the cases where $t_2$ or even $t_1 > s$ are easy modifications of the proof stated below.

\begin{rem}
Algorithm~\ref{alg} produces a vector
\begin{align*}
\tilde{\boldz} =  (\tilde{z}_1, \dots , \tilde{z}_s) = (z_1, \dots, z_{t_1},b^{w_{t_1+1}}z_{t_1+1}, \dots , b^{w_{t_2-1}}z_{t_2-1}, 0, \dots , 0)
\end{align*}
with the last $s - t_2 + 1$ components equal to 0. The elements in $b^{w_j}\znw{j}$ are evenly spread in $\zn$, where the notation $b^{w_j}\znw{j}$ means that each element of $\znw{j}$ is multiplied by $b^{w_j}$.
\end{rem}

As mentioned before we consider the squared worst-case error $\eworts$ defined above as a quality measure for the generating vector $\tilde{\boldz}$ produced by Algorithm~\ref{alg}. Thus we would like to find upper bounds for $\eworts$. A first bound can be directly deduced from the main result in \cite{DK}. To do so recall that we have defined the exclusion sets $\mathcal{E}_1, \dots, \mathcal{E}_s$ in Algorithm~\ref{alg} and consider the auxiliary exclusion sets $\tilde{\mathcal{E}_j} := \zn \setminus \znw{j} \cup \mathcal{E}_j$ for $j \in \{1, \dots , s \}$. Using these auxiliary exclusion sets in the projection-corrected algorithm constructed in \cite{DK} yields exactly the same vector as using the exclusion sets $\mathcal{E}_1, \dots, \mathcal{E}_s$ in Algorithm~\ref{alg} from the present paper. Replacing $z_j$ by $\tilde{z_j}$ does not hamper the proof in \cite{DK} and following it line by line we find the following upper bound on the squared worst-case error. Note that the coordinates $\tilde{z}_j$ of $\tilde{\boldz}$ do not necessarily belong to $\zn$. However, the formula for the squared worst-case error used in \cite{DK} is also true for arbitrary coordinates $\tilde{z}_j \in \{ 0, \dots, N-1\}$, see for example \cite{DKLP}.

\begin{thm}\label{thm:firstbound}
Let $s \in \NN$, $b \in \mathbbm{P}$, $m \in \NN$, $N = b^m$, $0 = w_1 \leq w_2 \leq \cdots $, and $\mathcal{Z}_N, \znw{j}, t_1$ and $t_2$ as above. Further let $\tilde{\boldz} = (\tilde{z}_1, \dots , \tilde{z}_s)$ be constructed by Algorithm~\ref{alg}, with exclusion sets $\mathcal{E}_j$ and auxiliary exclusion sets $\tilde{\mathcal{E}_j} := \zn \setminus \znw{j} \cup \mathcal{E}_j$ for $j \in \{1, \dots , s \}$. Then for all $1 \leq d \leq s$ we have
\begin{align*}
\ewortd \leq \left( \frac{1}{\phi(N)} \sum_{\mathfrak{u} \subseteq [d] }\gamma_{\mathfrak{u}}^{\lambda}(2 \zeta(\alpha \lambda))^{\left\lvert \mathfrak{u} \right\rvert } \prod_{j \in \mathfrak{u}} \frac{\phi(N)}{\phi(N) - \left\lvert \tilde{\mathcal{E}_j}\right\rvert } \right)^{\frac{1}{\lambda}}
\end{align*}
for any $\frac{1}{\alpha} < \lambda \leq 1$, where the product over the empty set is 1 by convention.
\end{thm}

\begin{rem}\label{remarkexclusion}
Note that the auxiliary exclusion sets used in this result can be considerably larger than the sets $\mathcal{E}_j$ actually used in Algorithm~\ref{alg} and thus the factor $\prod\limits_{j \in \mathfrak{u}} \frac{\phi(N)}{\phi(N) - \left\lvert \tilde{\mathcal{E}_j}\right\rvert }$ and with it the error bound are possibly unnecessarily large. Hence we want to adjust the proof so that we can use the exclusion sets $\mathcal{E}_j$, defined in Algorithm~\ref{alg} and thereby improve the error bound.
\end{rem}

Remark \ref{remarkexclusion} motivates the following

\begin{thm}\label{thm:betterbound}
Let $s \in \NN$, $b \in \mathbbm{P}$, $m \in \NN$, $N = b^m$, $0 = w_1 \leq w_2 \leq \cdots $, and $\mathcal{Z}_N, \znw{j}, t_1$ and $t_2$ as above. Further let $\tilde{\boldz} = (\tilde{z}_1, \dots , \tilde{z}_s)$ be constructed by Algorithm~\ref{alg}, with exclusion sets $\mathcal{E}_j$. Then for all $1 \leq d \leq s$ and $\frac{1}{\alpha} < \lambda \leq 1$ we have 
\begin{align*}
\ewortd \leq \left( \sum_{\mathfrak{u} \subseteq [d]}{ \frac{\gamma_{\mathfrak{u}}^{\lambda} (4 \zeta(\alpha \lambda))^{\left\lvert \mathfrak{u} \right\rvert} }{\phi(b^{\max\{ 0, m - \max_{j \in \mathfrak{u} }\{ w_j \} \}})} \prod_{j \in \mathfrak{u} }{\frac{\left\lvert \znw{j} \right\rvert}{\left\lvert \znw{j} \right\rvert - \left\lvert \mathcal{E}_j \right\rvert} } } \right)^{\frac{1}{\lambda}},
\end{align*}
where $\phi$ denotes the Euler totient function, and where we set $\max\emptyset = 0$.
\end{thm}

\begin{proof1}
The proof of \prettyref{thm:betterbound} is inspired by the proof in \cite{DK}. We use induction on $d$ to show the result. Recall our assumptions that $w_1 = 0$ and $z_1 = 1$, and that we consider product weights $\gamma_{\mathfrak{u}} = \prod_{j \in \mathfrak{u}}{\gamma_j}$. Then \eqref{eq:weight}, \eqref{eq:weight1dim}, and \eqref{eq:squaredWorst} together with Jensen's inequality, $\left( \sum_ka_k \right)^{\lambda} \leq \sum_ka_k^{\lambda}$ for non-negative $a_k$, yield
\begin{align*}
e^2(z_1)	&= \gamma_{\{1\}} \sum_{h \in \ZZ \setminus \{0 \}}{\left\lvert N h \right\rvert^{-\alpha}} = \frac{\gamma_{\{ 1 \}}}{N^{\alpha}} 2 \zeta(\alpha) \leq \left( \frac{\gamma_{\{ 1 \}}^{\lambda}}{b^m} 4 \zeta(\alpha \lambda) \right)^{\frac{1}{\lambda}} = \left( \frac{\gamma_{\{ 1 \}}^{\lambda} \left( 4 \zeta(\alpha \lambda) \right)^{\left\lvert \{ 1 \} \right\rvert} }{\phi(b^{\max\{ 0, m - \max_{j \in \{ 1 \}}w_j \}})} \right)^{\frac{1}{\lambda}}\\
					&= \left( \sum_{\mathfrak{u} \subseteq [1]}{ \frac{\gamma_{\mathfrak{u}}^{\lambda} (4 \zeta(\alpha \lambda))^{\left\lvert \mathfrak{u} \right\rvert} }{\phi(b^{\max\{ 0, m - \max_{j \in \mathfrak{u} }\{ w_j \} \}})} \prod_{j \in \mathfrak{u} }{\frac{\left\lvert \znw{j} \right\rvert}{\left\lvert \znw{j} \right\rvert - \left\lvert \mathcal{E}_j \right\rvert} } } \right)^{\frac{1}{\lambda}},
\end{align*}
as claimed.\newline

Now let $d \in \{ 1, \dots , s-1 \}$, and let $\tilde{z}_1, \dots , \tilde{z}_d$ be chosen with Algorithm~\ref{alg} and assume that
\begin{align*}
\ewortd \leq \left( \sum_{\mathfrak{u} \subseteq [d]}{ \frac{\gamma_{\mathfrak{u}}^{\lambda} (4 \zeta(\alpha \lambda))^{\left\lvert \mathfrak{u} \right\rvert} }{\phi(b^{\max\{ 0, m - \max_{j \in \mathfrak{u} }\{ w_j \} \}})} \prod_{j \in \mathfrak{u} }{\frac{\left\lvert \znw{j} \right\rvert}{\left\lvert \znw{j} \right\rvert - \left\lvert \mathcal{E}_j \right\rvert} } } \right)^{\frac{1}{\lambda}}
\end{align*}
holds for any $\lambda \in (\frac{1}{\alpha},1]$. We distinguish two cases, namely
\begin{enumerate}
	\item $d+1 \in \{ 2, \dots , t_1 \},$
	\item $d+1 \in \{ t_1 + 1, \dots , s \}.$
\end{enumerate}

Let us start with the first case where $d+1 \in \{ 2, \dots , t_1 \}$. As then $w_1 = \dots = w_{d+1} = 0$, we effectively consider the case of the projection-corrected CBC construction as in \cite{DK}. Thus we already know that
\begin{align*}
\ewortdd	&\leq \left( \frac{1}{\phi(N)} \sum_{\mathfrak{u} \subseteq [d+1]} \gamma_{\mathfrak{u}}^{\lambda} (2 \zeta(\alpha \lambda))^{\left\lvert \mathfrak{u} \right\rvert } \prod_{j \in \mathfrak{u}} \frac{\phi(N)}{\phi(N) - \left\lvert \mathcal{E}_j \right\rvert } \right)^{\frac{1}{\lambda}}\\
					&= \left( \sum_{\mathfrak{u} \subseteq [d+1]} \frac{\gamma_{\mathfrak{u}}^{\lambda} (4 \zeta(\alpha \lambda))^{\left\lvert \mathfrak{u} \right\rvert }}{\phi(b^{\max\{ 0, m - \max_{j \in \mathfrak{u} }\{ w_j \} \}})} \prod_{j \in \mathfrak{u}} \frac{\left\lvert \znw{j} \right\rvert}{\left\lvert \znw{j} \right\rvert - \left\lvert \mathcal{E}_j \right\rvert } \right)^{\frac{1}{\lambda}}
\end{align*}
for any $\lambda \in (\frac{1}{\alpha},1]$ and we are done with this case.

Next we deal with the second case where we have $d+1 \in \{ t_1 + 1, \dots , s \}$. Using \eqref{eq:squaredWorst} we easily obtain that for any $z \in \zn$
\begin{align*}
e^2(\tilde{z}_1 , \dots , \tilde{z}_d, z) = \ewortd + \theta_{N,d+1,\alpha, \bsgamma}(z),
\end{align*}
where
\begin{align*}
\theta_{N,d+1,\alpha, \bsgamma}(z) = \sum_{\substack{ \boldh \in \ZZ^{d+1} \\ h_{d+1} \neq 0 \\ \boldh \cdot (\tilde{z}_1, \dots , \tilde{z}_d, z) \equiv 0 \pmod{N}}}{r_{\alpha}(\bsgamma, \boldh)}.
\end{align*}

By setting $\beta_j = 1$ in \cite[Eq. (5)]{D}, we obtain
\begin{align}\label{eq:theta}
\theta_{N,d+1,\alpha, \bsgamma}(z) = 2 \gamma_{d+1} \zeta(\alpha) N^{-\alpha} (1 + \ewortd) + \gamma_{d+1} \kappa_{N, d+1, \alpha, \bsgamma}(z),
\end{align}
with
\begin{align}\label{eq:new}
\kappa_{N, d+1, \alpha, \bsgamma}(z) = \sum_{\substack{h_{d+1} \in \ZZ \\ N \nmid h_{d+1}}} \,\,\,\, \sum_{\substack{\boldh \in \ZZ^d \\ \boldh \cdot (\tilde{z}_1, \dots , \tilde{z}_d) \equiv -h_{d+1}z \pmod{N}}}{ \left\lvert  h_{d+1} \right\rvert^{-\alpha} r_{\alpha}(\bsgamma, \boldh)}.
\end{align}
Thus we have
\begin{align}\label{eq:pluszbeliebig}
\ewortdz	&= \ewortd + 2\gamma_{d+1}\zeta(\alpha) N^{-\alpha} (1 + \ewortd) + \gamma_{d+1} \kappa_{N,d+1,\alpha,\bsgamma}(z) \nonumber\\
					&= (1 + 2 \gamma_{d+1} \zeta(\alpha) N^{-\alpha}) \ewortd + 2 \gamma_{d+1} \zeta(\alpha) N^{-\alpha} + \gamma_{d+1}\kappa_{N,d+1,\alpha,\bsgamma}(z).
\end{align}
Recall that we want to show
\begin{align}\label{eq:zuzeigen}
\ewortdd \leq \left( \sum_{\mathfrak{u} \subseteq [d+1]}{ \frac{\gamma_{\mathfrak{u}}^{\lambda} (4 \zeta(\alpha \lambda))^{\left\lvert \mathfrak{u} \right\rvert} }{\phi(b^{\max\{ 0, m - \max_{j \in \mathfrak{u} }\{ w_j \} \}})} \prod_{j \in \mathfrak{u} }{\frac{\left\lvert \znw{j} \right\rvert}{\left\lvert \znw{j} \right\rvert - \left\lvert \mathcal{E}_j \right\rvert} } } \right)^{\frac{1}{\lambda}}
\end{align}
for any $\lambda \in (\frac{1}{\alpha},1]$.
 
Now choose $\lambda^* \in (\frac{1}{\alpha},1]$ such that the right hand side of \eqref{eq:zuzeigen} is minimized as a function of $\lambda$. 
Applying Jensen's inequality to \eqref{eq:pluszbeliebig} we obtain
\begin{align}\label{eq:error1}
(	&\ewortdz)^{\lambda^*}	\nonumber\\
	&\leq (1 + 2 \gamma_{d+1} \zeta(\alpha) N^{-\alpha})^{\lambda^*} (\ewortd)^{\lambda^*} + 2^{\lambda^*} \gamma_{d+1}^{\lambda^*}\zeta(\alpha)^{\lambda^*} N^{-\alpha\lambda^*} + \gamma_{d+1}^{\lambda^*}(\kappa_{N,d+1,\alpha,\bsgamma}(z))^{\lambda^*}\nonumber\\
	&\leq (1 + 2^{\lambda^*} \gamma_{d+1}^{\lambda^*} \zeta(\alpha\lambda^*) N^{-\alpha\lambda^*}) (\ewortd)^{\lambda^*} + 2^{\lambda^*} \gamma_{d+1}^{\lambda^*} \zeta(\alpha\lambda^*) N^{-\alpha\lambda^*} + \gamma_{d+1}^{\lambda^*}(\kappa_{N,d+1,\alpha,\bsgamma}(z))^{\lambda^*}.
\end{align}
Next we apply Jensen's inequality to \eqref{eq:new} and find
\begin{align*}
\frac{1}{\left\lvert \znw{d+1} \right\rvert} \sum_{l \in \znw{d+1}}(\kappa_{N,d+1,\alpha,\bsgamma}(l))^{\lambda^*} \leq \frac{1}{\left\lvert \znw{d+1} \right\rvert} \sum_{l \in \znw{d+1}}\kappa_{N,d+1,\alpha\lambda^*,\bsgamma^{\lambda^*}}(l) =: \overline{\kappa}_{N,d+1,\alpha\lambda^*,\bsgamma^{\lambda^*}}.
\end{align*}
In the following we use methods similar to \cite{DPW, DSWW}. Recall that Markov's inequality states that for a non-negative random variable $X$ and any real number $c \geq 1$ we have $\mathbbm{P}(X < c \mathbbm{E}(X)) > 1 - \frac{1}{c}$. We use the normalized counting measure $\mu$ on $\znw{d+1}$ as the probability measure and apply Markov's inequality as follows. For $c_{d+1} \geq 1$ let
\begin{align*}
G_{c_{d+1}} &:= \left\{ z \in \znw{d+1} \colon (\kappa_{N,d+1,\alpha,\bsgamma}(z))^{\lambda^*} \leq c_{d+1} \overline{\kappa}_{N,d+1,\alpha\lambda^*,\bsgamma^{\lambda^*}} \right\}\\
				&\supseteq \left\{ z \in \znw{d+1} \colon (\kappa_{N,d+1,\alpha,\bsgamma}(z))^{\lambda^*} \leq \frac{c_{d+1}}{\left\lvert \znw{d+1} \right\rvert} \sum_{l \in \znw{d+1}}{(\kappa_{N,d+1,\alpha,\bsgamma}(l))^{\lambda^*}} \right\} =: A_{c_{d+1}}.
\end{align*}
Then Markov's inequality yields
\begin{align*}
\mu(G_{c_{d+1}}) = \frac{\left\lvert G_{c_{d+1}} \right\rvert}{\left\lvert \znw{d+1} \right\rvert} \geq \mu(A_{c_{d+1}}) = \frac{\left\lvert A_{c_{d+1}} \right\rvert}{\left\lvert \znw{d+1} \right\rvert} > 1 - \frac{1}{c_{d+1}},
\end{align*}
that is, for any $c_{d+1} \geq 1$, exists a subset $G_{c_{d+1}} \subseteq \znw{d+1}$ of size bigger than $\left\lvert \znw{d+1} \right\rvert \left( 1 - \frac{1}{c_{d+1}} \right)$, such that
\begin{align*}
(\kappa_{N,d+1,\alpha,\bsgamma}(z))^{\lambda^*} \leq c_{d+1} \overline{\kappa}_{N,d+1,\alpha\lambda^*,\bsgamma^{\lambda^*}} \text{ for all } z \in G_{z_{d+1}}.
\end{align*}
By choosing $c_{d+1} \geq 1$ such that
\begin{align*}
\left\lvert \znw{d+1} \right\rvert \left( 1 - \frac{1}{c_{d+1}} \right) = \left\lvert \mathcal{E}_{d+1} \right\rvert,
\end{align*}
it is ensured that the set $G_{c_{d+1}} \setminus \mathcal{E}_{d+1}$ is not empty. Thus we have
\begin{align*}
c_{d+1} = \frac{\left\lvert \znw{d+1} \right\rvert}{\left\lvert \znw{d+1} \right\rvert - \left\lvert \mathcal{E}_{d+1} \right\rvert}.
\end{align*}
As $\tilde{z}_{d+1}$ is chosen by Algorithm~\ref{alg} such that the error $\ewortdd$ is minimal, we obtain together with \eqref{eq:error1}
\begin{align}\label{eq:error2}
\left(	\right.	& \left.\ewortdd \right)^{\lambda^*} \nonumber\\
								&\leq (1 + 2^{\lambda^*} \gamma_{d+1}^{\lambda^*} \zeta(\alpha\lambda^*) N^{-\alpha\lambda^*}) (\ewortd)^{\lambda^*} + 2^{\lambda^*} \gamma_{d+1}^{\lambda^*} \zeta(\alpha\lambda^*) N^{-\alpha\lambda^*} + \gamma_{d+1}^{\lambda^*}\overline{\kappa}_{N,d+1,\alpha\lambda^*,\bsgamma^{\lambda^*}} \nonumber\\
								&\leq (1 + c_{d+1} 2 \gamma_{d+1}^{\lambda^*} \zeta(\alpha\lambda^*) N^{-\alpha\lambda^*}) (\ewortd)^{\lambda^*} + c_{d+1} 2 \gamma_{d+1}^{\lambda^*} \zeta(\alpha\lambda^*) N^{-\alpha\lambda^*} + c_{d+1} \gamma_{d+1}^{\lambda^*}\overline{\kappa}_{N,d+1,\alpha\lambda^*,\bsgamma^{\lambda^*}}.
\end{align}
Using the induction assumption with $\lambda = \lambda^*$, we obtain
\begin{align}\label{eq:indassumpt}
 \left( \ewortd \right)^{\lambda^*} \leq \sum_{\mathfrak{u} \subseteq [d]}{ \frac{\gamma_{\mathfrak{u}}^{\lambda^*} (4 \zeta(\alpha \lambda^*))^{\left\lvert \mathfrak{u} \right\rvert} }{\phi(b^{\max\{ 0, m - \max_{j \in \mathfrak{u} }\{ w_j \} \}})} \prod_{j \in \mathfrak{u} }{\frac{\left\lvert \znw{j} \right\rvert}{\left\lvert \znw{j} \right\rvert - \left\lvert \mathcal{E}_j \right\rvert} } }.
\end{align}
Furthermore, from the proof of \cite[Lemma~5]{D}, we obtain
\begin{align}\label{eq:kappaquer}
\overline{\kappa}_{N,d+1,\alpha\lambda^*,\bsgamma^{\lambda^*}}	&\leq 2 \zeta(\alpha \lambda^*) (1 - N^{-\alpha \lambda^*}) \phi(N)^{-1} \sum_{\emptyset \neq \mathfrak{u} \subseteq [d]}{ \gamma_{\mathfrak{u}}^{\lambda^*} (2 \zeta(\alpha\lambda^*) )^{\left\lvert  \mathfrak{u} \right\rvert}} \nonumber\\
																																&\leq 2 \zeta(\alpha \lambda^*) (1 - N^{-\alpha \lambda^*}) \sum_{\emptyset \neq \mathfrak{u} \subseteq [d]}{ \frac{\gamma_{\mathfrak{u}}^{\lambda^*} (2 \zeta(\alpha\lambda^*) )^{\left\lvert  \mathfrak{u} \right\rvert}}{ \phi(b^{ \max \{ 0, m - \max_{j \in \mathfrak{u}} w_j \} }) } } \nonumber\\
																																&\leq 2 \zeta(\alpha \lambda^*) (1 - N^{-\alpha \lambda^*}) \sum_{\emptyset \neq \mathfrak{u} \subseteq [d]}{ \frac{\gamma_{\mathfrak{u}}^{\lambda^*} (4 \zeta(\alpha\lambda^*) )^{\left\lvert  \mathfrak{u} \right\rvert}}{ \phi(b^{ \max \{ 0, m - \max_{j \in \mathfrak{u}} w_j \} }) } \prod_{j \in \mathfrak{u}}{ \frac{\left\lvert \znw{j} \right\rvert}{\left\lvert \znw{j} \right\rvert - \left\lvert \mathcal{E}_j \right\rvert} } }.
\end{align}
Plugging \eqref{eq:indassumpt} and \eqref{eq:kappaquer} into \eqref{eq:error2} we have
\begin{align*}
\left( \right.	& \left. \ewortdd \right)^{\lambda^*} \leq\\
								&\leq (1 + c_{d+1} 2 \gamma_{d+1}^{\lambda^*} \zeta(\alpha\lambda^*) N^{-\alpha\lambda^*}) \sum_{\mathfrak{u} \subseteq [d]}{ \frac{\gamma_{\mathfrak{u}}^{\lambda^*} (4 \zeta(\alpha \lambda^*))^{\left\lvert \mathfrak{u} \right\rvert} }{\phi(b^{\max\{ 0, m - \max_{j \in \mathfrak{u} }w_j \}})} \prod_{j \in \mathfrak{u} }{\frac{\left\lvert \znw{j} \right\rvert}{\left\lvert \znw{j} \right\rvert - \left\lvert \mathcal{E}_j \right\rvert} } } \\
								&\quad + c_{d+1} 2 \gamma_{d+1}^{\lambda^*} \zeta(\alpha \lambda^*) N^{-\alpha \lambda^*} \\
								&\quad + c_{d+1} \gamma_{d+1}^{\lambda^*} 2 \zeta(\alpha \lambda^*) (1 - N^{-\alpha\lambda^*})\sum_{\emptyset \neq \mathfrak{u} \subseteq [d]}{ \frac{\gamma_{\mathfrak{u}}^{\lambda^*} (4 \zeta(\alpha\lambda^*) )^{\left\lvert  \mathfrak{u} \right\rvert}}{ \phi(b^{ \max \{ 0, m - \max_{j \in \mathfrak{u}} w_j \} }) } \prod_{j \in \mathfrak{u}}{ \frac{\left\lvert \znw{j} \right\rvert}{\left\lvert \znw{j} \right\rvert - \left\lvert \mathcal{E}_j \right\rvert} } }\\
								&= \sum_{\mathfrak{u} \subseteq [d]}{ \frac{\gamma_{\mathfrak{u}}^{\lambda^*} (4 \zeta(\alpha \lambda^*))^{\left\lvert \mathfrak{u} \right\rvert} }{\phi(b^{\max\{ 0, m - \max_{j \in \mathfrak{u} }w_j \}})} \prod_{j \in \mathfrak{u} }{\frac{\left\lvert \znw{j} \right\rvert}{\left\lvert \znw{j} \right\rvert - \left\lvert \mathcal{E}_j \right\rvert} } }\\
								&\quad + c_{d+1}\gamma_{d+1}^{\lambda^*} 2 \zeta(\alpha\lambda^*)N^{-\alpha \lambda^*}\frac{1}{\phi(N)}\\
								&\quad + c_{d+1}\gamma_{d+1}^{\lambda^*} 2 \zeta(\alpha\lambda^*)N^{-\alpha \lambda^*}\\
								&\quad + c_{d+1}\gamma_{d+1}^{\lambda^*} 2 \zeta(\alpha\lambda^*)\sum_{\emptyset \neq \mathfrak{u} \subseteq [d]}{ \frac{\gamma_{\mathfrak{u}}^{\lambda^*} (4 \zeta(\alpha\lambda^*) )^{\left\lvert  \mathfrak{u} \right\rvert}}{ \phi(b^{ \max \{ 0, m - \max_{j \in \mathfrak{u}} w_j \} }) } \prod_{j \in \mathfrak{u}}{ \frac{\left\lvert \znw{j} \right\rvert}{\left\lvert \znw{j} \right\rvert - \left\lvert \mathcal{E}_j \right\rvert} } }\\
								&\leq \sum_{\mathfrak{u} \subseteq [d]}{ \frac{\gamma_{\mathfrak{u}}^{\lambda^*} (4 \zeta(\alpha \lambda^*))^{\left\lvert \mathfrak{u} \right\rvert} }{\phi(b^{\max\{ 0, m - \max_{j \in \mathfrak{u} }w_j \}})} \prod_{j \in \mathfrak{u} }{\frac{\left\lvert \znw{j} \right\rvert}{\left\lvert \znw{j} \right\rvert - \left\lvert \mathcal{E}_j \right\rvert} } } \\
								&\quad + \frac{\left\lvert \znw{d+1} \right\rvert}{\left\lvert \znw{d+1} \right\rvert - \left\lvert \mathcal{E}_{d+1} \right\rvert}4 \gamma_{d+1}^{\lambda^*} \zeta(\alpha \lambda^*) N^{-\alpha \lambda^*} \\
								&\quad + \frac{\left\lvert \znw{d+1} \right\rvert}{\left\lvert \znw{d+1} \right\rvert - \left\lvert \mathcal{E}_{d+1} \right\rvert}4 \gamma_{d+1}^{\lambda^*} \zeta(\alpha \lambda^*) \sum_{\emptyset \neq \mathfrak{u} \subseteq [d]}{ \frac{\gamma_{\mathfrak{u}}^{\lambda^*} (4 \zeta(\alpha\lambda^*) )^{\left\lvert  \mathfrak{u} \right\rvert}}{ \phi(b^{ \max \{ 0, m - \max_{j \in \mathfrak{u}} w_j \} }) } \prod_{j \in \mathfrak{u}}{ \frac{\left\lvert \znw{j} \right\rvert}{\left\lvert \znw{j} \right\rvert - \left\lvert \mathcal{E}_j \right\rvert} } }\\
								&= \sum_{\mathfrak{u} \subseteq [d+1]}{ \frac{\gamma_{\mathfrak{u}}^{\lambda^*} (4 \zeta(\alpha \lambda))^{\left\lvert \mathfrak{u} \right\rvert} }{\phi(b^{\max\{ 0, m - \max_{j \in \mathfrak{u} }\{ w_j \} \}})} \prod_{j \in \mathfrak{u} }{\frac{\left\lvert \znw{j} \right\rvert}{\left\lvert \znw{j} \right\rvert - \left\lvert \mathcal{E}_j \right\rvert} } },
\end{align*}
as claimed. Thus the result holds for the special case of $\lambda^*$. As we have chosen $\lambda^*$ such that the right hand side of \eqref{eq:zuzeigen} is minimized the estimate is true for arbitrary $\lambda \in (\frac{1}{\alpha}, 1]$ as well.
\end{proof1}

\prettyref{thm:betterbound} enables us to combine the reduced with the projection-corrected CBC construction, while still achieving a small worst-case error. The reduced fast CBC construction, however, can be used equally well. In this case one has to perform the additional step of checking whether a component is in the respective exclusion set. Hence, in the process of choosing component $d$ one has to carry out $\left\lvert \mathcal{E}_d \right\rvert$ checks for exclusions at most, that is a total of at most
\begin{align*}
\left\lvert \mathcal{E}_2 \right\rvert + \dots + \left\lvert \mathcal{E}_s \right\rvert \leq (s - 1)\sum_{j = 2}^s\left\lvert \znw{j} \right\rvert \leq sN
\end{align*}
checks for the entire process of finding a generating vector. Hence the overall complexity of the reduced fast CBC algorithm is not increased.
\begin{rem}
As there is a connection between numerical integration using lattice point sets in the weighted Korobov space $\h$ and quasi-Monte Carlo integration using a randomly shifted or tent-transformed lattice point set in unanchored or anchored weighted Sobolev spaces, our results hold in these spaces as well. For detailed information see \cite{DKLP,DKS,H,HW,SKJ}.

It should also be possible to prove the results for general weights instead of just product weights.
\end{rem}

\end{document}